\documentclass[11pt]{amsart}
\usepackage{amsfonts,amsmath,amsthm,amscd,amssymb,latexsym,cite,verbatim,texdraw,floatflt,
caption2}
\RequirePackage{hyperref}
\usepackage[a4paper]{geometry}

\title{The normality of macrocubes and hyperballeans  }
\author{ Igor  Protasov, Ksenia Protasova}

\address{I.Protasov: Faculty of Computer Science and Cybernetics, Kyiv University,          Academic Glushkov pr. 4d, 03680 Kyiv, Ukraine}
\email{i.v.protasov@gmail.com}

\address{K.Protasova: Faculty of Computer Science and Cybernetics, Kyiv University,          Academic Glushkov pr. 4d, 03680 Kyiv, Ukraine}
\email{k.d.ushakova@gmail.com}

\begin{document}
\begin{abstract} For a  bornology $\mathcal B$ on a cardinal
$\kappa$, we  prove that
 the $\mathcal B$-macrocube is normal if and only if $\mathcal B$  has a  linearly ordered  base.  
As a corollary, we  get   that the   hyperballean of bounded subsets of an ultradiscrete ballean is not normal. These  answer Question 1 from \cite{b2}  and Question 14.4 from \cite{b1}.

\end{abstract}
\maketitle



\section{ Introduction  and preliminaries}

Given a set $X$, a family $\mathcal{E}$  of subsets of $X\times X$ is called a
{\it  coarse structure} on $X$ if

\begin{itemize}
\item{} each $E \in \mathcal{E}$  contains the diagonal $\bigtriangleup _{X}:=\{(x,x): x\in X\}$ of $X$;
\vspace{3 mm}

\item{}  if  $E$, $E^{\prime} \in \mathcal{E}$  then  $E \circ E^{\prime} \in \mathcal{E}$  and
$ E^{-1} \in \mathcal{E}$,    where  $E \circ E^{\prime} = \{  (x,y): \exists z\;\; ((x,z) \in E,  \ (z, y)\in E^{\prime})\}$,    $ E^{-1} = \{ (y,x):  (x,y) \in E \}$;
\vspace{3 mm}

\item{} if $E \in \mathcal{E}$ and  $\bigtriangleup_{X}\subseteq E^{\prime}\subseteq E$  then  $E^{\prime} \in \mathcal{E}$.
\end{itemize}

Elements $E\in\mathcal E$ of the coarse structure are called {\em entourages} on $X$

For $x\in X$  and $E\in \mathcal{E}$ the set $E[x]:= \{ y \in X: (x,y)\in\mathcal{E}\}$ is called the {\it ball of radius  $E$  centered at $x$}.
Since $E=\bigcup_{x\in X}\{x\}\times E[x]$, the entourage $E$ is uniquely determined by  the family of balls $\{ E[x]: x\in X\}$.
A subfamily ${\mathcal E} ^\prime \subseteq\mathcal E$ is called a {\em base} of the coarse structure $\mathcal E$ if each set $E\in\mathcal E$ is contained in some $E^\prime \in\mathcal E$.

The pair $(X, \mathcal{E})$  is called a {\it coarse space}  \cite{b10} or  a {\em ballean} \cite{b7}, \cite{b9}.

In this paper, all balleans under consideration are supposed to be
 {\it connected}: for any $x, y \in X$, there is $E\in \mathcal{E}$ such $y\in E[x]$.
A subset  $Y\subseteq  X$  is called {\it bounded} if $Y= E[x]$ for some $E\in \mathcal{E}$,
  and $x\in X$.
The family $\mathcal{B}_{X}$ of all bounded subsets of $X$  is a bornology on $X$.
We recall that a family $\mathcal{B}$  of subsets of a set $X$ is a {\it bornology}
if $\mathcal{B}$ contains the family $[X] ^{<\omega} $  of all finite subsets of $X$
 and $\mathcal{B}$  is closed   under finite unions and taking subsets. A bornology $\mathcal B$ on a set $X$ is called {\em unbounded} if $X\notin\mathcal B$.
A subfamily  $\mathcal B^{\prime}$ of $\mathcal B$ is called a base for $\mathcal B$ if, for each $B \in \mathcal B$, there exists $B^{\prime} \in \mathcal B^{\prime}$ such that $B\subseteq B^{\prime}$.

Each subset $Y\subseteq X$ defines a {\it subbalean}  $(Y, \mathcal{E}|_{Y})$  of $(X, \mathcal{E})$,
 where $\mathcal{E}|_{Y}= \{ E \cap (Y\times Y): E \in \mathcal{E}\}$.
A  subbalean $(Y, \mathcal{E}|_{Y})$  is called  {\it large} if there exists $E\in \mathcal{E}$
 such that $X= E[Y]$, where $E[Y]=\bigcup _{y\in Y} E[y]$.

Let $(X, \mathcal{E})$, $(X^{\prime}, \mathcal{E}^{\prime})$
 be  balleans. A mapping $f: X \to X^{\prime}$ is called
  {\it  macrouniform }  if for every $E\in \mathcal{E}$ there
  exists $E^{\prime}\in \mathcal{E}$  such that $f(E(x))\subseteq  E^{\prime}(f(x))$
    for each $x\in X$.
If $f$ is a bijection such that $f$  and $f ^{-1 }$ are macrouniform, then   $f  $  is called an {\it asymorphism}.
If  $(X, \mathcal{E})$ and  $(X^{\prime}, \mathcal{E}^{\prime})$  contains large  asymorphic  subballeans, then they are called {\it coarsely equivalent.}

Every  metric $d$ on a set  $X$ defines the coarse  structure
$\mathcal{E} _{d}$  on $X$  with  the base $\{\{(x,y): d(x, y)< n\} : n\in \mathbb{N}\}$.
A ballean $(X, \mathcal{E})$ is called {\it metrizable}  if there is a metric $d$ on such that $\mathcal{E} = \mathcal{E} _{d}$.

\vspace{5 mm}

{\bf Theorem 1.1. } ([9, Theorem 2.1.1])
{\it A ballean $(X,\mathcal{E})$ is metrizable  if and only if $\mathcal{E}$  has a countable base.}

\vspace{5 mm}

Let $(X,\mathcal{E})$  be a ballean. A subset $U \subseteq X$  is called an
 {\it asymptotic neighbourhood  } of a subset $Y\subseteq X$ if for every $E\in \mathcal{E}$ the set $E [Y] \setminus U$  is bounded.

Two subset $Y, Z$ of $X$  are called  {\it asymptotically disjoint  (separated)}  if for every $E\in \mathcal{E}$ the intersection 
  $E[Y] \cap E[Z] $  is bounded ($Y$ and $Z$ have disjoint asymptotic neighbourhoods).
We say that $Y,Z$  are {\it  linked} if $Y,Z$  are not asymptotically disjoint.

A ballean $(X, \mathcal{E})$  is called {\it normal}  \cite{b6} if any
two asymptotically disjoint subsets of $X$  are asymptotically separated. Every ballean $(X, \mathcal{E})$
 with linearly ordered base of $\mathcal{E}$  is  normal.
In particular,   every  metrizable ballean is normal, see [6, Proposition 1.1].

A function $f: X\to \mathbb{R}$  is called   {\it slowly oscillating} if 
 for any $E\in \mathcal{E}$  and $\varepsilon > 0$,  there exists a bounded  subset $B$ of $X$  such that $ diam \   f (E[x])<  \varepsilon $ for each $x\in X \setminus B$.   $ \ \  \Box $

\vspace{5 mm}

{\bf Theorem 1.2. } ([6, Theorem 2.1.])
{\it A ballean  $(X, \mathcal{E})$ is normal if and only if 
 for any two disjoint  asymptotically disjoint
subsets $Y, Z$ of $X$ there exists a slowly oscillating function  $f: X\to [0,1]$ such that
$f| _Y \equiv   0$  and $f| _Z  \equiv 1$.}
\vspace{5 mm}

{\bf Theorem 1.3. } ([1, Theorem 1.4 .])
{\it If the product $X \times Y$  of  two unbounded balleans $X, Y$ is normal then the  bornology  $B _ {X\times Y}$ has  a linearly  ordered base.  }

\section{ Macrocubes and hyperballeans}

Let     $\mathcal{B}$  be a bornology on an infinite cardinal  $\kappa$  such that $\kappa \notin \mathcal{B}$, 
$\{ 0,1 \}^{\mathcal{B}}=$
  $ \{ (x_{\alpha})_ {\alpha  < \kappa})\in  \{ 0,1 \}^\kappa :$ 
$\{ \alpha: x_\alpha =1\} \in \mathcal{B}. $
We take the family 
$$\{ \{(( x_\alpha)_ {\alpha  < \kappa},  \    (y_\alpha)_ {\alpha  < \kappa}) \in  \{ 0,1 \}^{\mathcal{B}  } \times
  \{ 0,1 \}^{\mathcal{B}}:  
x_{\alpha} = y_{\alpha}, \   \alpha\in \kappa \setminus B \} : B\in \mathcal{B}\}  $$
as a base  of the coarse structure on  $\{ 0,1 \}^{\mathcal{B}}$.
The obtained ballean is called the $\mathcal{B}$-macrocube \cite{b2}. 
If $\kappa = \omega$  and  $\mathcal{B}= [\omega]^{<\omega}$
then we get the well-known  Cantor macrocube whose coarse characterization was given in \cite{b3}.

Given a ballean $(X, \mathcal{E})$   the {\it hyperballean}   $(X, \mathcal{E})^\flat$   is a ballean on   $\mathcal{B}_{(X, \mathcal{E})}$ endowed  with  a coarse  structure  with the  base  $\{ E^\flat : E \in  \mathcal{E}\}, $    $E^\flat [Y]= \{Z: Y \subseteq E[Z]$,  $Z \subseteq E[Y]\}$, see \cite{b4},  \cite{b8}. 

Every bornology $\mathcal{B}$  on a cardinal $\kappa$  defines the {\it discrete ballean} 
$(\kappa,  \mathcal{E}_\mathcal{B} )$,  $ \mathcal{E}_\mathcal{B}$
has the base  
$\{  E_\mathcal{B}: B \in \mathcal{B} \}$, $E_\mathcal{B} [x]= B$
if 
 $x\in B$  and  
$E_\mathcal{B} [x]=\{x\}$
if  
$x\in \kappa \setminus B$.
If the set 
$\{\kappa \setminus B: B\in   \mathcal{B} \}$
 is an ultrafilter, then
$(\kappa,  \mathcal{E}_\mathcal{B} )$
 is called {\it ultradiscrete}.  
We denote by  
$[\kappa]^\mathcal{B}$
 the  hyperballean 
$(\kappa,  \mathcal{E}_\mathcal{B} )^\flat$.

\vspace{5 mm}

{\bf Proposition 2.1. }
{\it  Let $\kappa$  be a cardinal, $\mathcal{B} )$
be a  bornology on $\kappa$,  $K_0 = \{ A \in [\kappa]^\mathcal{B}  :  0\in A$.
Then the characteristic  function  $f: [\kappa]^\mathcal{B}  \longrightarrow \{0,1\}^\mathcal{B}$
is macro-uniform and the restriction
of $f$  to $K_0$ is an asymorphism between  $K_0$  and   $f(K_0)$. }

\vspace{5 mm}

{\bf Proposition 2.2. }
{\it   For any  cardinal  $\kappa$ and bornology  $\mathcal{B}$ on $\kappa$,  
the balleans
$\{0,1\}^\mathcal{B}$ and $[\kappa]^\mathcal{B}$  are not coarsely equivalent. 
\vspace{5 mm}

Proof.}
Following  \cite{b3}, we say that a ballean $(X, \mathcal{E}$ has asymptotically isolated balls if, for any bounded subset $B$ and  $E\in \mathcal{E}$  there exists $x\in X\setminus B$ such that  $E^\prime [x]\setminus E[x]= \emptyset$.
If  $(X^\prime , \mathcal{E} ^\prime)$ is  coarsely equivalent  to   $(X, \mathcal{E})$    then  $(X^\prime , \mathcal{E} ^\prime)$
has asymptotically isolated balls.

We observe
$\{0,1\}^\mathcal{B}$
 does not have asymptotically isolated balls, 
but the subballean
$[\kappa]_1 ^\mathcal{B} = \{ A\subset \kappa : |A|=1\}$
 of
$[\kappa] ^\mathcal{B} $
 has asymptotically isolated balls. Since
$[\kappa]_1 ^\mathcal{B} $
 and
$[\kappa] ^\mathcal{B}\setminus  [\kappa]_1 ^\mathcal{B}$
  are asymptotically disjoint, we see that
$[\kappa] ^\mathcal{B}$
 has asymptotically isolated balls. 
$ \ \  \Box $

\vspace{5 mm}

{\bf Proposition 2.3. }{\it
Let
$\mathcal{B} $
 be a bornology on a cardinal   $\kappa$  such that
$\{ \kappa \setminus B: B\in  \mathcal{B}\} $  is not an ultrafilter. Then
$[\kappa] ^\mathcal{B} $
 is normal if and only if
$\mathcal{B} $
  has a linearly ordered base. 

\vspace{5 mm}

Proof.}
If  $\mathcal{B} $ has a linearly ordered base then
$[\kappa] ^\mathcal{B} $
 is normal by  [6, Proposition 1.1].  

We assume that $[\kappa] ^\mathcal{B} $ is normal, partition $\kappa$  in two unbounded subsets $K_1$, $K_2$ and put 
$\mathcal{B}_1=\mathcal{B} | _{K_1}$, 
$\mathcal{B}_2=\mathcal{B} | _{K_2}$.
Then
$\{0,1\}^\mathcal{B}$
 is asymorphic to 
$\{0,1\}^{\mathcal{B}_1}\times   \{0,1\}^{\mathcal{B}_2}$.
By Theorem 1.3,  $\mathcal{B} $ has a linearly ordered base.
$ \ \  \Box $

\section{Some algebra on ultrafilters }

Let $G$ be an Abelian group.  We endow $G$ with the discrete topology and identify the Stone-$\check{C}$ech compactification $\beta G$ of  $G$ with the family of all ultrafilters on $G$. 

For $A\subseteq G$, we denote $\bar{A} =
 \{\mathcal{ U} \in \beta G: A\in \mathcal{ U}\}$.

Then the family $\{\bar{ A}: A\subseteq G \}$ forms a base for the topology on $\beta G$. 
Every mapping $f: G\longrightarrow [0,1]$ can be extended to the continuous  mapping
$ \  f^\beta : \beta G \longrightarrow [0,1]$.

Following \cite{b5}, we extend the addition  $+ $  on $G$ onto $\beta G$ by the following rule: for
$\mathcal{ U}$,  $\mathcal{ V}\in  \beta G$
 we take
$U\in \mathcal{ U}$
  and, for each $x\in U$, pick $V_x \in \mathcal{ U}$. 
Then $\bigcup_{x\in U} (x+V_x)\in \mathcal{ U}+\mathcal{V}$, and the family of all these subsets  forms a base of the ultrafilter 
$\mathcal{ U}+\mathcal{V}$.

\section{Results }

{\bf Theorem 4.1. }
{\it Let $ \mathcal{B}$ be a  bornology on a cardinal
$\kappa$.
The  $\mathcal B$-macrocube $X$ is normal if and only if $\mathcal B$ has a  linearly ordered base.
\vspace{5 mm}

Poof.} 
In light of Proposition 2.3,  it suffices to assume that $ \{ \kappa  \setminus B  : B\in \mathcal{B}\}$ is an ultrafilter   and  show that $X$  is not normal.

We consider   $X$  as a group with pointwise addition $mod \ 2$.
For $A \subseteq \kappa$ and  $\alpha  < \kappa$,  $\chi_A$  denotes the characteristic function of $A$,  $[0.\alpha]= \{\gamma < \kappa: \gamma \leq \alpha \}$.
Replacing $\kappa$ to some cardinal    $\kappa^\prime\leq  \kappa$,  we may suppose that $[0, \alpha]\in \mathcal B $ for  each 
$\alpha < \kappa$.

We consider two subsets $Y, Z$ of  $X$ defined by 
$$Y=\{ y_\alpha : \alpha < \kappa\},  \  y_\alpha  = \chi _{\{\alpha \}}  , $$
$$Z=\{ z_\alpha : \alpha < \kappa\},  \  z_\alpha  = \chi _{[0, \alpha]}  , $$
and show that $Y, Z$ are asymptotically disjoint. 
We take an arbitrary $B\in \mathcal B$,  denote by $H$ the subgroup 
of $X$ generated by
$\{ y_\alpha : \alpha  \in  B \} $
 and take the minimal  
$\gamma < \kappa$
such that $\gamma \in \kappa\setminus B$. 
Then
$(y_\beta  + H)\cap Z = \emptyset $
  as  soon as $\beta \in \kappa\setminus B$, $\beta > \gamma$.

We suppose that $X$ is normal and use Theorem 1.2  to  choose a slowly oscillating function $f: X\longrightarrow [0,1]$ 
such that 
$f|_{Z}  \equiv  1$, 
$f|_{Y}  \equiv  0$,
and denote by 
$\mathcal{U} $  and $\mathcal{V} $ 
ultrafilters on $X$ with the bases 
$$\{\{  y_\alpha : \alpha\in  \kappa \setminus B\} : B\in \mathcal B\}, \  \{\{  z_\alpha : \alpha\in  \kappa \setminus B\} : B\in \mathcal B\}. $$

Let $f^\beta (\mathcal{V}  + \mathcal{U}) = r  $,   $r\in [0.1]. $ 
We take an arbitrary 
$W\in \mathcal{V}  + \mathcal{U} $ 
and pick  $ C \in \mathcal B$
and
 $ D _\alpha \in  \mathcal B  $,
$\alpha \in  C  $
such that
$$\bigcup \{ z_\alpha +y_\beta:  \alpha \in \kappa\setminus    C,  \ 
 \beta \in \kappa\setminus  D_\alpha \}\subseteq W. $$
Then we construct inductively two mappings
$\psi: \kappa \setminus  C\longrightarrow \kappa,  $ 
$\phi: \kappa \setminus  C\longrightarrow \kappa $ 
such that  
$\psi(\alpha) \in \kappa \setminus  D_\alpha $, $\phi(\alpha)\in \kappa\setminus  D_\alpha  $ 
and 
$\psi (\kappa \setminus  C) \cap \phi(\kappa\setminus  C) =\emptyset$. 
Since 
$\{  \kappa \setminus B : B\in \mathcal B\}$
is an ultrafilter, either 
$\psi(\kappa \setminus  C)\in \mathcal B$ or
$\phi(\kappa \setminus  C)\in \mathcal B$.
We assume that 
$\psi(\kappa \setminus   C)\in \mathcal B$
 and denote by $H$ the subgroup of $X$ generated by 
$\{ y_\alpha: \alpha \in \psi(\kappa \setminus  C)\}.$

Then 
$z_\alpha \in  W+H$
 for each 
$\alpha \in  \kappa \setminus  C$.
Hence, $W$ and $Z$ are linked. 
Since $f$ is slowly oscillating, we conclude that $r=1$.

On the other hand, for 
$\alpha\in \kappa \setminus  C$,
 $ z_\alpha +  \{y_\beta:   
 \beta \in \kappa\setminus  D_\alpha \}\subseteq W. $
It follows that $W$ and $Y$ are  linked. 
Since $f$ is slowly oscillating, we get $r=0$ contradicting above paragraph. $\ \Box$

\vspace{5 mm}

{\bf Corollary 4.2. }
{\it  Let $\mathcal B$  be a bornology on a cardinal $\kappa$ such that 
$\{  \kappa \setminus B : B\in \mathcal B\}$
 is an ultrafilter. 
Then the hyperballean
$[\kappa ]^{\mathcal B}$
  is not normal. 

\vspace{5 mm}

Poof.}  Since a subballean of a normal ballean is normal, to apply Theorem 4.1, we use Proposition 2.1.  $ \ \Box$
\vspace{5 mm}

Theorem 4.1 answers Question 1 from \cite{b2}, 
Corollary 4.2 answers  Question 14.4 from \cite{b1}.
\vspace{5 mm}

For  a bornology  $\mathcal B$ on a cardinal  $\kappa$, the subballean of  all characteristic functions of finite subsets of 
$\kappa$ of the 
ballean $\{ 0,1\}^\mathcal B$
is called the {\it finitary $\mathcal B$-macrocube. }
It follows from the proof of  Theorem 4.1 that the finitary $\mathcal B$-macrocube on $\omega$ is not normal provided that $\{ \omega\setminus B : B\in \mathcal B\}$ is ultrafilter.

\vspace{5 mm}

{\bf Question 4.3. }
{\it  Let $\mathcal B$  be a bornology on a cardinal $\kappa$ such that 
$\{  \kappa \setminus B : B\in \mathcal B\}$
 is an ultrafilter. 
Is it true that the finitary $\mathcal B$-macrocube is not normal?}

\vspace{5 mm}

For $n\in \omega$, $n>0$ and a bornology   $\mathcal{B}$  on $\kappa$, we denote  
 $[\kappa]^\mathcal{B}_n  = \{ A\in [\kappa]^\mathcal{B} : |A|\leq n\}$,  
 $\{ 0,1 \}^\mathcal{B}_n  = \{  (x_\alpha)_{\alpha<\kappa}: | \{ \alpha : x_  \alpha  =1\} |  \leq n \}$.
If  $\{ \kappa  \setminus B: B \in \mathcal{B} \} $ is an ultrafilter then 
 $[\kappa]^\mathcal{B}_n $ is normal  [1, Theorem 1.13]. 

A ballean $(X, \mathcal{E})$ is called {\it ultranormal} if any two unbounded subsets of $X$ are linked.

\vspace{3 mm}

{\bf Theorem 4.4. }
{\it Let $ \mathcal{B}$ be  a  bornology on $\kappa$  such that 
$\{ \kappa \setminus B: B \in \mathcal{B}\}$  is an ultrafilter. 
Then the ballean 
 $\{ 0,1 \}^\mathcal{B}_n  $ is ultranormal. 
\vspace{5 mm}

Proof.}  We consider  $\{ 0,1 \}^\mathcal{B}_1  $ as a subballean of    $\{ 0,1 \}^\mathcal{B}_n  $.
Since  
 $\{ \kappa  \setminus B: B \in \mathcal{B} \} $  is an ultrafilter, it suffices to show that $\mathcal{A}$, 
$\{ 0,1 \}^\mathcal{B}_1  $ are linked for any unbounded subset $\mathcal{A}$ of $\{ 0,1 \}^\mathcal{B}_n  $.

Applying $n$-times Lemma 7.1 from \cite{b1}, we pick $P\subseteq \kappa$ such that $ \kappa \setminus P\in \mathcal{B}$ 
and  $ \mathcal{A}^\prime = \{ A \in  \mathcal{A}:  |A\cap P|=1\}$ is unbounded.
 We put 
$H= \{ (x_\alpha) _ {\alpha <\kappa} : x_\alpha =0$ for each $\alpha \in P\}$.
Then $(x+ H) \cap  \mathcal{A}^\prime \neq \emptyset$ for each $x\in P$. 
Hence, $ \mathcal{A}^\prime$  and 
 $\{ 0,1 \}^\mathcal{B} _1  $ are linked.  $ \ \  \Box $



\begin{thebibliography}{10}

\bibitem{b1} T.~Banakh, I.~Protasov, {\em The normality and bounded growth of balleans},\newline {\tt https://arxiv.org/abs/1810.07979}.


\bibitem{b2} T.~Banakh, I.~Protasov, {\em Constructing  balleans}, J. Math. Sciences,  {\bf 241} (2019), 16-26.


\bibitem{b3}{ T. Banakh, I. Zarichnyi,}  {\it Characterizing  the Cantor bi-cube in asymptotic categories,} Groups Geom. Dyn. {\bf 5}: 4 (2011), 691-728.


\bibitem{b4} D. Dikranjan,  I.  Protasov, K.  Protasova,  N.  Zava,  {\em Balleans, hyperballeans and ideals}, Appl. Gen. Topology
{\bf 2}: 4 (2019), 431-447.



\bibitem{b5}{  N. Hindman, D. Strauss, }{\it Algebra in the Stone-$\check{C}$ech Compactification}, de Gructer, Berlin, New York, 1998.

\bibitem{b6} { I. Protasov, }{\it Normal ball structures},  Math. Stud. {\bf 20}    (2003), 3-16.


\bibitem{b7}{ I. Protasov, T. Banakh, }{\it Ball Structures and Colorings of Groups and Graphs},  Math. Stud. Monogr. Ser., Vol. 11, VNTL, Lviv, 2003.


\bibitem{b8}{ I. Protasov,  K. Protasova,} {\it  On  hyperballeans of bounded geometry    },  Europ. J. Math.
 {\bf 4}  (2018),  1515-1520.

\bibitem{b9}{ I. Protasov, M. Zarichnyi,} {\em General Asymptology},   Math. Stud. Monogr. Ser., Vol. 12, VNTL, Lviv, 2007,  pp. 219.


\bibitem{b10}    {J. Roe,} {\em Lectures on Coarse Geometry}, Univ. Lecture Ser., vol. 31, American Mathematical Society, Providence RI, 2003.



\end{thebibliography}
\end{document}